\documentclass{amsart}

\usepackage{amsmath}
\usepackage{amssymb}
\usepackage{amsthm}
\usepackage{eucal}

\theoremstyle{plain}
\newtheorem{thm}{Theorem}[section]
\newtheorem{prp}[thm]{Proposition}

\newtheorem{fct}[thm]{Fact}
\newtheorem{cor}[thm]{Corollary}

\theoremstyle{definition}
\newtheorem{dfn}[thm]{Definition}

\theoremstyle{remark}
\newtheorem{rmk}[thm]{Remark}

\newtheorem{qst}[thm]{Question}

\newtheorem{clm}{Claim}[thm]

\newtheorem{subclm}[clm]{Subclaim}

\newtheorem{clm_lem}[clm]{Lemma}

\newtheorem{tens_equ}[clm]{\m\bigotimes}

\theoremstyle{definition}

\newenvironment{prf}{\begin{proof}}{\end{proof}}

\newcommand{\m}{\ensuremath}

\renewcommand{\l}{\m{\ell}}
\newcommand{\lam}{\m{\lambda}}
\newcommand{\om}{\m{\omega}}
\newcommand{\ka}{\m{\kappa}}

\newcommand{\ph}{\m{\varphi}}

\newcommand{\et}{\m{\eta}}
\renewcommand{\ni}{\m{\nu}}
\newcommand{\thet}{\m{\vartheta}}

\newcommand{\then}{\m{\Rightarrow}}
\renewcommand{\implies}{\m{\Longrightarrow}}
\renewcommand{\iff}{\m{\Longleftrightarrow}}
\newcommand{\sat}{\m{\models}}
\newcommand{\set}[1]{\ensuremath{\{#1\}}}
\newcommand{\seq}[1]{\ensuremath{\langle#1\rangle}}

\newcommand{\card}[1]{\ensuremath{|#1|}}
\newcommand{\rest}{\m{\upharpoonright}}
\newcommand{\cont}{\m{\frown}}

\newcommand{\initeq}{\m{\unlhd}}

\newsavebox{\indbinbox}
\savebox{\indbinbox}
{\begin{picture}(0,0)
    \newlength{\gnu}
    \settowidth{\gnu}{$\smile$}
    \setlength{\unitlength}{.5\gnu}
    \put(-1,-.675){$\smile$}
    \put(-.275,.1){$|$}
  \end{picture}}

\newsavebox{\indbinAbox}
\savebox{\indbinAbox}
{\begin{picture}(0,0)
    \settowidth{\gnu}{$\smile$}
    \setlength{\unitlength}{.5\gnu}
    \put(-1,-.675){$\smile$}
    \put(-.275,.1){$|$}
    \put(-.95,.2){$\to$}
  \end{picture}}

\DeclareMathOperator{\tp}{tp}
\DeclareMathOperator{\rk}{rk}

\DeclareMathOperator{\len}{len}

\DeclareMathOperator{\tand}{\; and \;} \DeclareMathOperator{\tfor}{\; for \;}
 
\DeclareMathOperator{\tif}{\; if \;}

\newcommand{\FC}{\m{\mathfrak{C}} \;}

\newcommand{\CL}{\m{\mathcal{L}} \;}

\newcommand{\CB}{\m{\mathcal{B} \;}}

\newcommand{\elementary}{\m{\prec}}

\newcommand{\initial}{\lhd}
\newcommand{\mr}{\begin{itemize}}
\newcommand{\ermn}{\end{itemize}}
\newcommand{\roster}{\begin{itemize}}
\newcommand{\eendroster}{\end{itemize}}
\newcommand{\itemitem}{\item}
\newcommand{\ub}{\underline}
\newcommand{\definition}{\begin{dfn}}
\newcommand{\eenddefinition}{\end{dfn}}
\newcommand{\proclaim}{\begin{prp}}
\newcommand{\eendproclaim}{\end{prp}}
\newcommand{\scite}{\cite}
\newcommand{\wilog}{\text{wlog}}
\newcommand{\dbcu}{\bigcup}
\newcommand{\Cal}{\mathcal}
\newcommand{\dsize}{}
\newcommand{\ifff}{\iff}
\newcommand{\initialeq}{\unlhd}
\newcommand{\iso}{\approx}
\newcommand{\deq}{=:}

\newcommand{\x}{\m{\bar x}}
\newcommand{\y}{\m {\bar y}}
\renewcommand{\c}{\m {\bar c}}
\renewcommand{\a}[1]{\mbox{\m{\bar a_{#1}}}}

\newcommand{\aet}{\mbox{\m{\bar a_\et}}}

\newcommand{\treeq}[1]{\mbox{\m{^{{#1} \ge }2}}}

\newcommand{\phix}[1]{\mbox{\ensuremath{\varphi(\x, #1)}}}
\newcommand{\phixa}[1]{\mbox{\ensuremath{\varphi(\x, \a{#1})}}}

\newcommand{\phixasup}[2]{\mbox{\ensuremath{\varphi^{#2}(\x, \a{#1})}}}

\newcommand{\aettreeq}[1]{\mbox{\m{\set{\a{\et} : \et \in
\treeq{#1}}}}}

\newcommand{\phibranchalt}[2]{\mbox{\m{\set{\phixasup{#1 \rest
i}{#1(i)} : i < #2}}}}

\newcommand{\dcont}[2]{\mbox{\ensuremath{#1 \frown #2}}}

\newcommand{\dcontz}[1]{\mbox{\ensuremath{#1 \frown \langle 0 \rangle}}}
\newcommand{\dconto}[1]{\mbox{\ensuremath{#1 \frown \langle 1 \rangle}}}

\newcommand{\qy}{\mbox{\ensuremath{q(\bar y)}}}
\newcommand{\px}{\mbox{\ensuremath{p(\bar x)}}}

\newcommand{\nonex}[1]{\mbox{\m{{\neg (\exists
\bar{x})(\varphi(\bar{x},\bar{y})
\wedge\varphi(\bar{x},#1))}}}}

\newcommand{\cupph}[2]{\mbox{\m{#1\cup\{\varphi(\bar{x},
#2)\}}}}
\newcommand{\pcupph}[1]{\mbox{\m{\cupph{\px}{#1}}}}
\newcommand{\cupex}[2]{\mbox{\m{#1 \cup\{\neg (\exists
\bar{x})(\varphi(\bar{x},\bar{y})
\wedge\varphi(\bar{x},#2))\}}}}
\newcommand{\qcupex}[1]{\mbox{\m{\cupex{\qy}{#1}}}}

\newcommand{\rank}[2]{\mbox{\m{\rk^1_\ph(#1,#2)}}}
\newcommand{\rkpq}{\mbox{\m{\rk^1_\ph(\px,\qy)}}}

\newcommand{\sop}[1]{\mbox{\m{SOP_{#1}}}}
\newcommand{\sopa}{\mbox{\m{SOP_{1}}}}

\newcommand{\sopt}{\mbox{\m{SOP'_{1}}}}

\def\vp{\varphi}
\def\w{\omega}

\def\uhr{\upharpoonright}

\def\lor{\langle 1\rangle}
\def\gc{\mathfrak{C}}
\def\and{\&}

\title{More on SOP$_{1}$ and SOP$_{2}$}
\author{Saharon Shelah \and
    Alex Usvyatsov}
\address{
Mathematics Department\\
Hebrew University of Jerusalem\\
91904 Givat Ram, Israel\\
}

\begin{document}

\begin{abstract}
This paper continues \cite{DjSh692}. We present a rank
function for NSOP$_{1}$ theories and give an example of a
theory which is NSOP$_{1}$ but not simple. We also
investigate the connection between maximality in the
ordering $\initial^*$ among complete first order theories
and the (N)SOP${}_2$ property. We complete the proof
started in
\cite{DjSh692} of the fact that
$\initial^*$-maximality implies SOP${}_2$ and get weaker
results in the other direction. The paper provides a step
toward the classification of unstable theories without the
strict order property.
\end{abstract}

\maketitle

\section{Introduction and preliminaries}
This paper continues \scite{DjSh692} and investigates
theories that have or do not have the order properties
\sop{1} and \sop{2}. These properties were
defined in \scite{DjSh692} in order to find more division
lines lying between the tree property (non-simplicity) and
\sop{3}, the first dividing line in
Shelah's hierarchy of finite approximations of the strict
order property. We remind the definitions:

Let $T$ be a complete first order theory, \FC - the monster
model of $T$ (a $\ka^*$ - saturated and homogeneous model
for $\ka^*$ big enough).

\begin{dfn}
\begin{enumerate}
\item
    Let $n \ge 3$. We say \phix{\y} (with $\len(x) = \len(y)$)
    exemplifies the strong order property of order $n$ (\sop{n})
    in $T$ if it
    defines on \FC a graph with infinite indiscernible
    chains with no cycles of length $n$.
\item
    We say \phix{\y} (with $\len(x) = \len(y)$)
    exemplifies the \emph{strict order property} in $T$ if it
    defines on \FC a partial order with infinite indiscernible
    chains.
\end{enumerate}
\end{dfn}

\begin{fct}
    For a theory $T$,
    strict order property \implies \sop{n+1} \implies
    \sop{n} for all $n \ge 3$.
\end{fct}
\begin{prf}
    The first implication is trivial, for the other one see
    \cite{Sh500}, claim (2.6).
\end{prf}

We also remind an equivalent definition of \sop{3}:
\begin{fct}
$T$ has \sop{3} if and only if there is an indiscernible
sequence $\langle
\bar{a}_i:\,i<\omega\rangle$ and formulae
$\varphi(\bar{x},\bar{y})$, $\psi(\bar{x},\bar{y})$ such
that
\begin{itemize}
\item[(a)] $\{\varphi(\bar{x},\bar{y}), \psi(\bar{x},\bar{y})\}$ is contradictory,
\item[(b)] for some sequence $\langle \bar{b}_j:\,j<\omega\rangle$ we have
\[
i\le j\implies \models\varphi [\bar{b}_j,\bar{a}_i]\mbox{
and } i> j\implies \models\psi [\bar{b}_j,\bar{a}_i].
\]
\item[(c)] for $i<j$, the set $\{\varphi(\bar{x},\bar{a}_j),
\psi(\bar{x},\bar{a}_i)\}$ is contradictory.
\end{itemize}
\end{fct}

\begin{prf}
    Easy, or see \scite{Sh500}, claim (2.20).
\end{prf}

Now we recall the definitions of \sop{1}, \sop{2} and
related properties:

\begin{dfn}\label{dfn:nsop21}
\begin{enumerate}
\item $T$ has $SOP{}_2$ if there is a formula $\varphi(\bar{x},
\bar{y})$ which exemplifies this property in
$\mathfrak{C}$, and this means:

There are $\bar{a}_\eta\in {\mathfrak C}$ for $\eta\in
{}^{\omega>}2$ such that
\begin{itemize}
\item[(a)] For every
$\eta\in {}^{\omega}2$, the set $\{\varphi(\bar{x},
\bar{a}_{\eta\rest\l}):\,l<\omega\}$ is consistent.
\item[(b)] If $\eta,\nu\in {}^{\omega>}2$ are incomparable,
$\{\varphi(\bar{x},\bar{a}_\eta),
\varphi(\bar{x},\bar{a}_\nu)\}$ is inconsistent.
\end{itemize}

\item $T$ has $SOP{}_1$ if there is a formula $\varphi(\bar{x},
\bar{y})$ which exemplifies this in ${\mathfrak C}$, which means:

There are $\bar{a}_\eta\in {\mathfrak C}$, for $\eta\in
{}^{\omega>}2$ such that:

\begin{itemize}
\item[(a)] for $\rho\in {}^\omega 2$ the set
$\{\varphi(\bar{x},\bar{a}_{\rho\rest n}):
\,n<\omega\}$ is consistent.

\item[(b)] if $\nu\frown\langle 0\rangle\initialeq \eta\in {}^{\omega>}2$,
\underline{then} $\{\varphi(\bar{x},\bar{a}_\eta), \varphi(\bar{x},
\bar{a}_{\nu\frown\langle 1\rangle})\}$ is inconsistent.

\end{itemize}

\item  $NSOP{}_2$ and $NSOP{}_1$ are the negations of $SOP{}_2$ and
$SOP{}_1$ respectively.

\item  $T$ has $SOP'{}_1$ if there is a formula $\varphi(\bar{x},
\bar{y})$ which exemplifies this property in
$\mathfrak{C}$, and this means:

there are $\langle\bar{a}_\eta:\,\eta\in {}^{\omega>} 2
\rangle$ in ${\mathfrak C}_T$ such that
\begin{itemize}
\item[(a)] $\{\varphi(\bar{x},\bar{a}_{\eta\rest n})^{\eta(n)}:\,
n<\omega\}$ is consistent for every $\eta\in {}^\omega 2$,
where we use the notation
\[
\varphi^l=\left\{
\begin{array}{cc}
\varphi &\mbox{ if }l=1,\\
\neg\varphi &\mbox{ if }l=0
\end{array}
\right.
\]
for $l<2$.
\item[(b)] If $\nu\frown\langle 0\rangle\initialeq\eta\in {}^{\omega>}2$,
\underline{then} $\{\varphi(\bar{x},\bar{a}_\eta), \varphi(\bar{x},
\bar{a}_\nu)\}$ is inconsistent.
\end{itemize}

\item $T$ has $SOP''{}_2$ if there is a formula $\varphi(\bar{x},
\bar{y})$ which exemplifies this property in
$\mathfrak{C}$, and this means:

there is a sequence
\[
\left
\langle
\bar{a}_{\bar{\eta}}:\,\bar{\eta}=\langle\eta_0,\ldots\eta_{n-1}\rangle,
\eta_0\initial\eta_1\initial\ldots\initial\eta_{n-1}\in {}^{\lambda>}
2 \mbox{ and }\lg(\eta_i)\mbox{ successor}
\right\rangle
\]
such that
\begin{itemize}
\item[($a$)] for each $\eta\in {}^\lambda 2$, the set
\begin{equation*}
\left\{
\begin{array}{ll}
\ph(x,\bar{a}_{\bar{\eta}}):\,&\bar{\eta}=
\langle\eta\rest(\alpha_0+1),
\eta\rest(\alpha_1+1),\ldots \eta
\rest(\alpha_{n-1}+1)\rangle\\
&\mbox{ and }
\alpha_0<\alpha_1<\ldots\alpha_{n-1}<\lambda\\
\end{array}
\right\}
\end{equation*}
is consistent

\item[($b$)] for every large enough $m$,
if $h$ is a 1-to-1 function from $^{n \ge} m$ into
${}^{\lambda
>}2$ preserving $\eta \initial \nu$ and $\nu
\perp
\nu$ (incomparability) \ub{then} $\{\varphi(\bar x,\bar a_{\bar
\nu}):\text{for some } \eta \in {}^n m \text{ we have } \bar \nu =
\langle h(\eta \restriction \ell):\ell \le n \rangle\}$ is
inconsistent.

\end{itemize}

\end{enumerate}
\end{dfn}

\begin{fct}\label{fct:impl}
\begin{enumerate}
\item
    For a theory $T$, \sop{3} \implies \sop{2} \implies
    \sop{1}
\item
    $T$ has \sop{1} if and only if it has \sopt
\end{enumerate}
\end{fct}
\begin{prf}
    See \cite{DjSh692}.
\end{prf}

It is still not known whether the implications in
~\ref{fct:impl}, (1)  are strict, but for now we
investigate each one of these order properties on its own.

In the second section we expand our knowledge on \sop{1}.
We present a rank function measuring type-definable
``squares'', i.e. pairs of types of the form $(\px,\qy)$
and show the rank is finite for every such a pair if and
only if $T$ does not have \sopt (if and only if $T$ does
not have \sop{1}). In fact, if one calls a tree of
parameters  $\set{\bar{a}_\eta : \eta\in {}^{\omega>}2}$
showing that \phix{\y} exemplifies \sopt in \FC (as in the
definition of \sopt) a $\ph-\sopt \,tree$, the rank
measures exactly the maximal depth of a tree like this that
can be built in \FC. We also show a small application of
the rank.

It is easy to see (see \cite{DjSh692}) that if \phix\y \;
exemplifies \sop{1} in \FC then it also exemplifies the
tree property, so $T$ has \sop{1} \implies $T$ is not
simple. We show that the implication is proper, i.e. find
an example of a theory $T$ which is not simple, but is
$NSOP_1$. This theory which we call $T^*_{\rm feq}$, was
first defined in \cite{Sh457}, and is used in
\cite{Sh500} as an example of an $NSOP_3$ non-simple
theory. Here we use a slightly different definition of the
same theory, as given in \cite{DjSh692}.

\begin{dfn}
\begin{enumerate}
\item
    $T_{\rm feq}$ is the following theory
in the language $\{Q, P, E, R, F\}$
\begin{itemize}
\item[(a)]
    Predicates $P$ and $Q$ are unary and disjoint, and
$(\forall x)\,[P(x)\vee Q(x)]$,
\item[(b)]
    $E$ is an equivalence relation on $Q$,
\item[(c)]
    $R$ is a binary relation on $Q\times P$ such that
\[
[x\,R\,z\,\,\&\,\,y\,R\,z\,\,\&\,\,x\,E\,y]\implies x=y.
\]
{\scriptsize{(so $R$ picks for each $z\in Q$ (at most one)
representative of any $E$-equivalence class).}}
\item[(d)]
    $F$ is a (partial) binary function from $Q\times P$ to $Q$, which
satisfies
\[F(x,z)\in Q\,\,\&\,\, \,\left(F(x,z)\right)\,R\,z \,\,\&\,\,
x\, E\,\left( F(x,z)\right).
\]
{\scriptsize{(so for $x\in Q$ and $z\in P$, the function
$F$ picks the representative of the $E$-equivalence class
of $x$ which is in the relation $R$ with $z$).}}
\end{itemize}

\item
     $T^\ast_{\rm feq}$ is the model completion
of $T_{\rm feq}$, (so a complete theory with infinite
models, in which $F$ is a full function).

\end{enumerate}
\end{dfn}

If the reader thinks about the definition above, he'll find
out that $T^*_{\rm feq}$ is just the model completion of
the theory of infinitely many (independent) parametrised
equivalence relations. The reader can also compare between
the definition of $T^*_{\rm feq}$ here and in \cite{Sh457}.
As we already mentioned, it was shown in \cite{Sh500} this
theory does not have \sop{3} (but is not simple). Here we
prove an (a priori) stronger result: $T^*_{\rm feq}$ does
not have \sop{1}.

In the third section we deal with
$\initial^*_{\lambda}$-maximality (see the beginning of the
section for definitions). For a theory $T$, to be
$\initial^*_{\lambda}$-maximal means to be complicated. In
a way, it means that it is hard to make its models
$\lambda$-saturated. In \cite{Sh500} it was stated that
\sop{3} implies $\initial^*_{\lambda}$-maximality; here we
fill the missing details of the proof, showing explicitly
that the model completion of the theory of trees is
$\initial^*_{\lambda}$-maximal for every regular \lam \,
big enough.

We are interested in this paper, though, in \sop{2} more
than \sop{3}. In \cite{DjSh692} it was shown that a
property similar to $\initial^*_{\lambda}$-maximality
(which also follows from $\initial^*_{\lambda}$-maximality
for some \lam \, under certain set theoretic conditions)
implies $SOP''_2$, and one of the questions asked there is
what is the connection between $SOP''_2$ and the \sop{n}
hierarchy. Of course, it would be natural to try to connect
between $SOP''_2$ and \sop{2}, and indeed we prove here
these two properties are equivalent for a theory $T$ (not
necessarily for a formula).

So we can conclude
\sop{3} \implies $\initial^*_{\lambda}$-maximality \implies
\sop{2}. Unfortunately, we don't know much about the other
directions of the above implications.

In \cite{DjSh692} two notions of ``tree indiscernibility''
were defined. We recall the definitions:

\begin{dfn}\label{dfn:indisc}
\begin{enumerate}
\item Given an ordinal $\alpha$ and sequences $\bar{\eta}_l=
\langle \eta^l_0,\eta^l_1,\ldots,\eta_{n_l}^l\rangle$ for $l=0,1$
of members of ${}^{\alpha>} 2$, we say that
$\bar{\eta}_0\iso_1\bar{\eta}_1$ iff
\begin{itemize}
\item[(a)] $n_0=n_1$,
\item[(b)] the truth values of
\[
\eta^l_{k_3}\initialeq \eta^l_{k_1}\cap\eta^l_{k_2},\quad
\eta^l_{k_1}\cap\eta^l_{k_2}\initial \eta^l_{k_3}, \quad
(\eta^l_{k_1}\cap\eta^l_{k_2})\frown\langle
0\rangle\initialeq \eta^l_{k_3},
\]
for $k_1,k_2, k_3\le n_0$, do not depend on $l$.
\end{itemize}

\item We say that the sequence $\langle\bar{a}_\eta:\,
\eta\in {}^{\alpha>}2\rangle $ of $\mathfrak C$ (for an
ordinal $\alpha$) are {\em 1-fully binary tree
indiscernible (1-fbti)} iff whenever
$\bar{\eta}_0\iso_1\bar{\eta}_1$ are sequences of elements
of ${}^{\alpha>}2$, \underline{then}
\[
\bar{a}_{\bar{\eta}_0} \deq \bar{a}_{\eta^0_0}\frown\ldots \frown\bar{a}_{\eta^0_{n_0}}
\]
and the similarly defined $\bar{a}_{\bar{\eta}_1}$, realize
the same type in $\mathfrak C$.

\item We replace 1 by 2 in the above definition
iff $(\eta^l_{k_1}\cap\eta^l_{k_2})
\frown \langle 0\rangle\initialeq\eta^l_{k_3}$ is omitted from clause (b)
above.
\end{enumerate}
\end{dfn}

We will need the following fact proved in \cite{DjSh692},
(2.11):

\begin{fct}\label{fct:thinning} If
$t\in \{1,2\}$ and $\langle \bar{b}_\eta:\,\eta\in
{}^{\omega>} 2 \rangle$ are given, and $\delta\ge \omega$,
\underline{then} we can find $\langle
\bar{a}_\eta:\,\eta\in {}^{\delta>}2\rangle$ such that
\begin{itemize}
\item[(a)] $\langle \bar{a}_\eta:\,\eta\in {}^{\delta>}2\rangle$ is $t$-fbti,
\item[(b)] if $\bar{\eta}=\langle \eta_m:\,m<n\rangle$, where each $\eta_m
\in {}^{\delta>}2$ is given, and $\Delta$ is a finite set
of formulae of $T$, \underline{then} we can find $\nu_m\in
{}^{\omega>}2\,(m<n)$ such that with $\bar{\nu}\deq
\langle \nu_m:\,m<n\rangle$, we have $\bar{\nu}\iso_t\bar{\eta}$
and the sequences $\bar{a}_{\bar{\eta}}$ and
$\bar{b}_{\bar{\nu}}$, realise the same $\Delta$-types.
\end{itemize}
\end{fct}

\section{More on \sop{1}}
We work with a complete first order theory $T$, let \FC be
its ``monster'' model (saturated in some very big $\ka^*$).
Let $\CL = \CL(T)$ (the language of $T$). Every formula we
mention in this section is an \CL-formula, maybe with
parameters from \FC.

First, we would like to make sure that we indeed are 
developing a new theory here. As every simple theory is $NSOP_1$,
it is very important to ask whether the other direction is also true
(if so, we would find ourselves in a well-developed context, for which
almost all the theorems proven here are either known or easy). 
But the answer is negative:

\begin{thm}
$T^*_{feq}$ does not have $SOP_1$.
\end{thm}
\begin{prf}
Suppose there exists $\vp(\bar x, \bar y)$ with $\ell g (\bar x) =
n,\;  \ell g(\bar y) = m$, and $\langle \bar a_\eta: \eta \in^{\w
>} 2\rangle$ in ${}^m\gc$ which exemplify $SOP_1$ in $\gc$ ($\gc$ is
the monster model of $T^*_{feq})$. Without loss of generality,
(by ref$\{$fct:thinning$\}$)$\langle\bar a_\eta: \eta \in {}^{\w
> } 2\rangle$ if 1-full tree indiscernible.  Also, by elimination
of quantifiers, we may assume that $\varphi(\bar x, \bar y)$ is
quantifier free.  As the only function symbol in the language is
$F$ and $F^\gc$ has the property $F^\gc(F^\gc(x, z), y) = F^\gc
(x, y)$ for all $z$, we will also assume wlog that $\bar x$ and
$\bar y$ in $\varphi (\bar x, \bar y)$ are closed under $F$ and
$\varphi (\bar x, \bar y)$ gives the full diagram of $\bar x
\frown \bar y$.  We shall regard $\bar x $ as $\langle x^0,\dots,
x^{n-1}\rangle, \; \bar y $ as $\langle y^0,\dots,
y^{m-1}\rangle,\;  \bar a_\eta$ as $\langle a^0_\eta,\dots,
a^{m-1}_\eta\rangle$.

By the definition of $SOP_1$, there exist $\bar e = \langle
e^0,\dots,e^{n-1}\rangle, \; \bar d = \langle d^0,\dots,
d^{n-1}\rangle$ in ${}^n\gc$ s.t.
\begin{equation*}
\gc\models  \vp(\bar e, \bar a_{\langle\; \rangle})\land \vp(\bar
e, \bar a_{\langle 0 \rangle}) \land \vp(\bar e, \bar a_{\langle
00 \rangle})
\end{equation*}
and
\begin{equation*}
\gc\models\vp(\bar d, \bar a_{\langle\; \rangle})\land \vp(\bar d,
\bar a_{\langle 1 \rangle})
\end{equation*}

Denote $\eta= \langle 00\rangle$.  Let $B = \gc \uhr\bar a_\eta
\frown\bar a_{\langle 1 \rangle}$.  By our assumptions, there
exists a model $N_0$ whose universe is $\bar x \frown \bar
a_\eta$, extending $\gc\uhr\bar a_\eta$, whose basic diagram is
$\vp(\bar x, \bar a_\eta)$.  Similarly, there exists a model
$N_1$ with universe $\bar x \frown \bar a_{\lor})$ and basic
diagram $\vp(\bar x, \bar a_{\lor})$.
 We shall
amalgamate $B, N_0$ and $N_1$ into a model of $T_{feq}, N$.  This
will immediately give a contradiction:  first, extend $N$ to
$N^*\models T^*_{feq}$, then amalgamate $N^*$ and $\gc$ over $B$
into some $\gc^+\models T^*_{feq}$.  By model completeness of
$T^*_{feq}$, $\gc \prec \gc^+$, but $\gc^+\models \exists \bar
x(\vp(\bar x, \bar a_\eta)\land \vp (\bar x, \bar a_{\lor}))$,
which is a contradiction to the definition of $SOP_1$.

It is left, therefore, to show that we can define on $|N_0|\cup
|N_1|$ a structure which will be a model of $T_{feq}$, extending
$B$.

We define $N$ as follows:
\begin{equation*}
|N|=|N_1|\cup|N_2|,\quad P^N=P^{N_1} \cup P^{N_2},\qquad  Q^N =
Q^{N_1} \cup Q^{N_2}.
\end{equation*}

Note that the diagram of $\bar x$ in $N_0$ is the same as the
diagram of $\bar x$ in $N_1$ (both implied by $\vp(\bar x, \bar
y)$, and the diagrams of $\bar a_\eta, \bar a_{\lor}$ in $N_i$ are
the same as in $\gc$, hence the same as in $B$.  Therefore, $P^N$
and $Q^N$ are well defined and give a partition of $|N|$.  Also,
so far $N$ extends $B$ (as a structure).

Considering $E$ and $R$, we define
\begin{align*}
R^N &= R^{N_1} \cup R^{N_2} \cup R^B\\
E^N &= E^{N_1}\cup E^{N_2} \cup E^B
\end{align*}
Once we have proven the following lemmas, we will be able to define
$F^N$ in a natural way, and in fact will be done.

\begin{clm_lem}\label{lem:1}  
$E^N$ is an equivalence relation on $Q^N$,
extending $E^B$.
\end{clm_lem}

\begin{clm_lem}\label{lem:2}
$R^N$ is a two-place relation on $N$, \ \ $R^N\subseteq P^N\times
Q^N$, satisfying:

for every $y \in P^N$ and every equivalence class $C$ of $E^N$,
there exists a unique $z \in C$ such that $(y, z)\in R^N$.
\end{clm_lem}

\bigskip
\noindent{\it Proof of \ref{lem:1}}. The only nonobvious thing is
transitivity.  We check two main cases, all the rest are either
similar or trivial.
\begin{enumerate}
\item[(1)]  Assume $x^iE^Na^j_\eta, \; \; x^iE^Na^k_{\lor}$ for some $i,
j, k$.  We want to show $a^j_\eta E^N a^k_{\lor}$.  It is enough
to see $a^j_\eta E^{\gc} a^k_{\lor}$.  We will write $E$ instead
of $E^{\gc}$. \newline $N\models x^iEa^j_\eta\Rightarrow N_0
\models x^i E a^j_\eta \Rightarrow \vp(\bar x, \bar y) \vdash
x^iEy^j$.  Similarly, $\vp (\bar x, \bar y) \vdash x^i Ey^k$, and
we get (by the choice of $\bar e, \bar d \in {}^n\gc$) 

$e^i Ea^j_\eta,\; e^i Ea^j_{\langle \; \rangle}, \; e^iEa^k_\eta, \;
e^iEa^k_{\langle \; \rangle}, d^i Ea^j_{\lor},\; d^i
Ea^j_{\langle\; \rangle}\;, d^iEa^k_{\lor}, \; d^iEa^k_{\langle
\; \rangle}$. 

Now it is easy to see that all the above elements are $E$-equivalent
in $\gc$, in particular $a^j_\eta$ and $a^k_{\lor}$, as required. 

\item[(2)]  Assume $x^iE^Na^n_\eta, \; a^k_{\lor} E^Na^n_\eta$, and
we show $x^iE^Na^k_{\lor}$, i.e. $\vp(\bar x, \bar y) \vdash x^i
Ey^k$.  As $\vp (\bar d, \bar a_{\lor})$ holds in $\gc$ and
$\vp(\bar x, \bar y)$ gives a full diagram, it will be enough to
see $d^iEa^k_{\lor}$.
\end{enumerate}

We know that $\vp(\bar x, \bar y) \vdash x^i Ey^j$ therefore
$e^iEa^j_\eta, \; e^iEa^j_{\langle \;\rangle},\; d^iEa^j_{\lor},
\; d^iEa^j_{\langle\;\rangle}$.  In particular, $d^i Ea^j_\eta$,
but, by our assumption, $a^j_\eta Ea^k_{\lor}$, so we are done.

\hfill $\square_1$

\bigskip
\noindent{\it Proof of \ref{lem:2}}. Like in the previous lemma, the
only nontrivial thing to prove is the last part, and we will deal
with two main cases.

\begin{enumerate}
\item $N\models (a^i_\eta R\, a^j_{\lor}) \land (a^i_\eta Rx^k)
\land (x^kEa^j_{\lor})$. We aim to show $N\models
(x^k=a^j_{\lor})$. We know: 
\begin{itemize}
\item[$(*)_1$] $\gc \models a^i_\eta R\, a^j_{\lor}$
\item[$(\ast)_2$] $N_0 \models a^i_\eta R x^k$, therefore $\vp
(\bar x, \bar y) \vdash y^iRx^k$
\item[$(\ast)_3$] $N_1 \models x^k Ea^j_{\lor}$, therefore $\vp
(\bar x, \bar y) \vdash x^k Ey^j$.
\end{itemize}

So we can conclude:

\begin{itemize}
\item[$(\ast)_2$] $\Rightarrow a^i_{\langle\; \rangle} Re^k,\; \;
\; a^i_{\langle\;\rangle} Rd^k$
\item[$(\ast)_3$] $\Rightarrow e^kEa^j_{\langle\; \rangle}, \; \;
\; d^kEa^j_{\langle\; \rangle} \Rightarrow e^kEd^k$.
\end{itemize}

As the above two relations hold in $\gc$, which is a model of
$T_{feq}$, we get $\gc\models e^k=d^k$.  Denote $e^\ast =
e^k=d^k$.

\begin{itemize}
\item[$(\ast)_1$] $\Rightarrow a^i_\eta R\, a^j_{\lor}$

\item[$(\ast)_2$] $\Rightarrow a^i_\eta R e^\ast$

\item[$(\ast)_1$] $\Rightarrow e^\ast E a^j_{\lor}$
\end{itemize}
Together (once again, $\gc \models T_{feq}$) we get $e^\ast =
a^j_{\lor}$, therefore $\vp(\bar x, \bar y) \vdash x^k=y^j$, so
$N_1\models x^k = a^j_{\lor}$, and we are done.
\item $N\models (x^i R\, a^j_{\lor}) \land (x^i R\, a^k_\eta) \land (a^k_\eta E a^j_{\lor})$
 and we aim to show $N\models (a^k_\eta= a^j_{\lor})$.

 We know:
 \begin{itemize}
 \item[$(\ast)_1$] $N_1 \models x^i R\, a^j_{\lor}$, so $\vp(\bar x,
 \bar y) \vdash x^i Ry^j$
 \item[$(\ast)_2$] $N_0 \models x^i R\, a^k_\eta$, so $\vp (\bar x,
 \bar y) \vdash x^i Ry^k$
 \item[$(\ast)_3$] $\gc \models a^k_\eta Ea^j_{\lor}$
 \end{itemize}
 \end{enumerate}
Note that by indiscernibility of $\langle \bar a_r: r \in {}^{w>}
2\rangle$ and $(\ast)_3$ we get $a^k_{\langle 0 \rangle}
Ea^j_{\lor}$, therefore $a^k_{\langle 0 \rangle} Ea^k_\eta$.  Now,
by $(\ast)_2$, $ e^iR\, a^k_\eta \; \;  \& \; \;
e^iR\,a^k_{\langle 0 \rangle}$. Therefore, by $\gc \models
T_{feq}$, $a^k_{\langle 0 \rangle} = a^k_\eta$. Now by
indiscernibility
\begin{equation*}
a^k_{\langle 0 \rangle}= a^k_{\langle \;\rangle}, \quad
a^k_{\lor}= a^k_{\langle \;\rangle}
\end{equation*}
So we get that all of the above are equal (and in fact $a^k_{r_1}
= a^k_{r_2}$ for all $r_1, r_2 \in {}^{w>} 2)$.

Now:
\begin{itemize}
\item[$(\ast)_1$] $\Rightarrow d^iR \,a^j_{\lor}$
\item[$(\ast)_2$] $\Rightarrow d^iR\, a^k_{\lor} \Rightarrow
d^iR\, a^k_\eta$ (as $a^k_{\lor} = a^k_\eta)$
\item[$(\ast)_3$] $\Rightarrow a^k_\eta Ea^j_{\lor}$.
\end{itemize}
By $\gc \models T_{feq}$, we conclude $a^k_\eta = a^k_{\lor}$,
which finishes the proof of the lemma, and therefore the proof of
the theorem. \hfill $\square_2$ 

\end{prf}

Our next goal is to show that there is a rank function
closely connected with being (N)SOP${}_1$. Let \ph{\x,\y}
be a formula.

\begin{dfn}\label{rank} Given (partial) types $p(\bar{x}),
q(\bar{y})$. By induction on $n<\omega$ we define when
\[
\mbox{rk}{}^1_{\varphi(\bar{x},\bar{y})}(p(\bar{x}), q(\bar{y}))\ge n:
\]

If \underline{$n=0$}, this happens if both $p(\bar{x}),
q(\bar{y})$ are consistent

For \underline{$n+1$}, the rank is $\ge n+1$ if for some
\c \, \sat \qy, both
\[
\mbox{rk}{}^1_{\varphi(\bar{x},\bar{y})}(p(\bar{x})\cup\{\varphi(\bar{x},
\bar{c})\}, q(\bar{y}))\ge n
\]
and
\[
\mbox{rk}{}^1_{\varphi(\bar{x},\bar{y})}(p(\bar{x}), q(\bar{y})\cup\{\neg (\exists
\bar{x})(\varphi(\bar{x},\bar{y}) \wedge\varphi(\bar{x},\bar{c}))\})\ge n.
\]

We say rk${}^1_{\varphi(\bar{x},\bar{y})}(p(\bar{x}), q(\bar{y}))=\infty$ iff
rk${}^1_{\varphi(\bar{x},\bar{y})}(p(\bar{x}), q(\bar{y}))\ge n$ for all $n$.

We say the rank is $-1$ if it is not bigger or equal to
$0$.
\end{dfn}

\begin{rmk}\label{rmk:rank}
(1) The statement rk${}^1_{\varphi(\bar{x}, \bar{y})}(\theta_1(\bar{x};
\bar{a}),\theta_2(\bar{x}; \bar{b}))\ge n$ is a first order formula with parameters
$\bar{a},\bar{b}$.

{\noindent (2)} We can continue to define when
rk${}^1_{\varphi(\bar{x},\bar{y})}(p(\bar{x}), q(\bar{y}))\ge \alpha$ for any
ordinal $\alpha$, but by the compactness theorem and part (1) it follows that
rk${}^1_{\varphi(\bar{x},\bar{y})}(p, q)\ge \alpha$ for some $\alpha\ge\omega$ iff
rk${}^1_{\varphi(\bar{x},\bar{y})}(p, q)\ge \omega$ iff
rk${}^1_{\varphi(\bar{x},\bar{y})}(p, q) = \infty$.

{\noindent (3)} (Monotonicity) If $p'\vdash p''$ and $q'\vdash q''$,
\underline{then} rk${}^1_{\varphi(\bar{x},\bar{y})}(p', q')\le
$rk${}^1_{\varphi(\bar{x},\bar{y})}(p'', q'')$.


{\noindent (4)} (Finite Character) If rk${}^1_{\varphi(\bar{x},\bar{y})}
(p(\bar{x}), q(\bar{y}))=n$, \underline{then} for some finite $p_0(\bar{x})\subseteq
p(\bar{x})$ and $q_0(\bar{y})\subseteq q(\bar{y})$ we have
rk${}^1_{\varphi(\bar{x},\bar{y})} (p_0(\bar{x}), q_0(\bar{y}))=n$.


{\noindent (5)} If $p'\equiv p''$, and $q'\equiv q''$, \underline{then}
rk${}^1_\varphi(p',q')=$rk${}^1_\varphi(p'',q'')$.


\end{rmk}

We aim to show that \rkpq is finite for every \px, \qy (or,
equivalently, \rank{\x = \x}{\y = \y} is finite) if and
only if \ph(\x,\y) does not exemplify \sopt in $T$. For
this purpose we shall need another definition and several
easy claims.

\begin{dfn}\label{def:soptree}
    Given (partial) types \px and \qy,
    we say that \set{\aet: \et \in \treeq{n}} is a
    \emph{\ph-\sopt tree for \px and \qy \, (of depth $n$)} if
    \begin{itemize}
        \item[(a)] $\px \cup \phibranchalt{\et}{n}$ is
        consistent for every $\et \in ^n2$.
        \item[(b)] $\aet \models \qy$ for all $\et \in
        \treeq{n}$
        \item[(c)] If \et, \ni \, are in \treeq{n}
        satisfying \dcontz{\et} \initeq \ni, then the set
        \set{\phixa{\et}, \phixa{\ni}} is inconsistent.
    \end{itemize}
\end{dfn}

\begin{prp}\label{prp:specialtree}
    Suppose \aettreeq{n} is a \ph-\sopt tree for \px and \qy \, of depth
    $n$, and denote $A^0 = \set{\aet : \seq{0} \initeq
    \et}$, $A^1 = \set{\aet : \seq{1} \initeq
    \et}$. Then
    \begin{enumerate}
    \item
        $A^1$ is a \ph-\sopt tree for \pcupph{\a{\seq{}}}
        and \qy
    \item
        $A^0$ is a \ph-\sopt tree for \px and
        \qcupex{\a{\seq{}}}.
    \end{enumerate}
\end{prp}
\begin{prf}

    The clauses (a) and (c) of the definition easily
    hold both for $A^1 \tand A^0$,
    so we should only check (b), which is also
    obvious for $A^1$. Therefore, we're left to show that
    for every $\et \in A^0$, $\aet \models
    \nonex{\a{\seq{}}}$, and this is clear by clause (c) of
    the definition (\aettreeq{n} is a \ph-\sopt \, tree, and
    \mbox{\dcont{\seq{}}{0} \initeq \et}).
\end{prf}

Now we show the connection between the rank and \sopt
trees.

\begin{prp}\label{prp:ranktree}
    $\rkpq \ge n \iff$ there exists a \ph-\sopt tree for
    \px and \qy of depth $n$.
\end{prp}
\begin{prf}
    Both directions are proved by induction on $n$. The
    case $n=0$ is obvious. For $n=m+1$, the right-to-left
    direction follows immediately by the induction
    hypothesis and ~\ref{prp:specialtree}. So we will
    elaborate more only about the other direction, although
    it is also straightforward.

    Suppose $n=m+1$ and $\rkpq \ge n$. By the definition
    of the rank and the induction hypothesis,
    for some \c \, \sat \qy, there are
    \begin{enumerate}
    \item
        a \ph-\sopt
        tree $A^1 = \set{\a{\et}^1: \et \in \treeq{m}}$
        for \pcupph{\c} and \qy
    \item
        a \ph-\sopt
        tree $A^0= \set{\a{\et}^0: \et \in \treeq{m}}$
        for \px and \qcupex{\c}
    \end{enumerate}
    (both of depth $m$). We define a tree \aettreeq{n} by
    \[
    \begin{array}{ll}
        & \a{\seq{}} = \c
            \\
        & \a{\seq{\l} \cont \et} =
        \a{\et}^{\l} \tfor
        \l \in \set{0,1}
    \end{array}
    \]
    which is as required, i.e. a \ph-\sopt tree for \px and \qy. Why?
    \begin{itemize}
    \item[(a)]
        of the definition obviously holds by (1)
        above.
    \item[(b)] holds as \c \sat \qy.
    \item[(c)] obviously holds by (2) above.
    \end{itemize}
\end{prf}

The following remark is obvious:
\begin{rmk}
    \phix{\y} exemplifies \sopt in $T$ \iff there exists a
    \ph-\sopt tree for $\x=\x \tand \y=\y$ of any depth.
\end{rmk}

So we can conclude the following

\begin{thm}
    A formula \phix{\y} does not exemplify \sopt in $T$ \iff
    \mbox{$\rank{\x=\x}{\y=\y} < \om$} \iff
    \mbox{$\rkpq < \om$} for every two (partial) types \px and
    \qy. Moreover, \rank{\x=\x}{\y=\y} is exactly the maximal
    depth of a \ph-\sopt tree that can be built in \FC.

\end{thm}

\begin{cor}
    $T$ does not have \sopa \iff $T$ does not have \sopt
    \iff \rank{\x=\x}{\y=\y} is finite for every formula
    \phix{\y}.
\end{cor}

Now we show an application of the rank.

\begin{thm}\label{thm:rnkapl} Suppose that $T$ satisfies NSOP${}_1$. Assume that
\begin{itemize}
\item[(a)] $M_1\elementary M_2\elementary \mathfrak C$.
\item[(b)] $p$ is a (not necessarily complete)
type over $M_2$, containing the formula
$\varphi(\bar{x},\bar{b}^\ast)$ for some $\bar b^\ast\in
M^2\setminus M^1$.
\end{itemize}
\underline{Then} for some
finite $q' \subseteq \tp(\bar b^\ast / M_1)$
at least one of the following holds:
\begin{itemize}
\item[(i)] If $\bar{b}\in M_1$ realises $q'(\bar{y})$ \underline{then}
$\varphi(\bar{x},\bar{b})\notin p$, or
\item[(ii)] If $\bar{b}\in M_1$ realises $q'(\bar{y})$ then
$\{\varphi(\bar{x}, \bar{b}),
\varphi(\bar{x},\bar{b}^\ast)\}$ is consistent.
\end{itemize}

In fact, all we need to assume for this Claim is that
$\varphi(\bar{x},\bar{y})$ does not exemplify that $T$ is
SOP${}_1$.
\end{thm}

\begin{prf} Denote $q = {\rm tp}(b^\ast/M_1)$.
As $T$ is NSOP${}_1$,we have that
rk${}^1_{\varphi(\bar{x},\bar{y})}(p\rest M_1, q)=n^\ast<
\omega$ (certainly $n^\ast \ge 0$).
By the finite character of the rank, we have that for some
finite $p_0\subseteq p\rest M_1$ and $q_0\subseteq q$,
\[
\mbox{rk}{}^1_{\varphi(\bar{x},\bar{y})}(p_0, q_0)=n^\ast.
\]
Hence for no $\bar{c} \models q_0(\y)$ do we have that both
$\mbox{rk}{}^1_{\varphi(\bar{x},\bar{y})}(p_0\cup\{\varphi(\bar{x},\bar{c})\},
q_0)\ge n^\ast$ and
$\mbox{rk}{}^1_{\varphi(\bar{x},\bar{y})}(p_0,
q_0\cup\{(\neg\exists\bar{x})[\varphi(\bar{x},\bar{y})
\wedge\varphi(\bar{x},\bar{c})]\})\ge n^\ast$.
In particular, this holds for $\bar{c}=\bar{b}^\ast$
(remember that $\bar{b}^\ast
\models q$ and therefore certainly $\bar{b}^\ast \models q_0$).
So
\begin{tens_equ}
\label{tens:bstar}
\[
    \mbox{\;\; If } \mbox{rk}{}^1_{\varphi(\bar{x},\bar{y})}
    (p_0\cup\{\varphi(\bar{x},\bar{b}^\ast)\},q_0)\ge
    n^\ast, \mbox{ then }
    \mbox{rk}{}^1_{\varphi(\bar{x},\bar{y})}(p_0,
    q_0\cup\{(\neg\exists\bar{x})[\varphi(\bar{x},\bar{y})
    \wedge\varphi(\bar{x},\bar{b}^\ast)]\}) < n^\ast.
\]
\end{tens_equ}

 By Remark \ref{rmk:rank}(1), there is a finite
$q'\subseteq q$ such that
\begin{tens_equ}
\label{tens:qtag}
\[
\bar{b}\mbox{ realises }q'\implies
\mbox{rk}{}^1_{\varphi(\bar{x},\bar{y})}(p_0\cup\{\varphi(\bar{x},\bar{b})\},
q_0)=\mbox{rk}{}^1_{\varphi(\bar{x},\bar{y})}
(p_0\cup\{\varphi(\bar{x},\bar{b^\ast })\}, q_0).
\]
\end{tens_equ}
We aim to show that $q'$ is as required.

\underline{Case 1}. $\mbox{rk}{}^1_{\varphi(\bar{x},\bar{y})}
(p_0\cup\{\varphi(\bar{x},\bar{b^\ast})\}, q_0)=n< n^\ast$.

We note that the possibility (i) holds.


Namely, suppose $\bar{b}$ realises $q'$,then
$\mbox{rk}{}^1_{\varphi(\bar{x},\bar{y})}(p_0\cup\{\varphi(\bar{x},\bar{b})\},
q_0)=n< n^\ast$, so if $\varphi(\bar{x},\bar{b})\in p$, we
obtain a contradiction with monotonicity of the rank.

\underline{Case 2}. $\mbox{rk}{}^1_{\varphi(\bar{x},\bar{y})}
(p_0\cup\{\varphi(\bar{x},\bar{b^\ast})\}, q_0)=n^\ast$.

We shall show that (ii) holds.


Suppose otherwise, so let $\bar{b}\in M_1$ realise $q'$ and
$\{\varphi(\bar{x},\bar{b}),
\varphi(\bar{x},\bar{b^\ast})\}$ is contradictory. By
~\ref{tens:qtag},
\[\mbox{rk}{}^1_{\varphi(\bar{x},\bar{y})}
(p_0\cup\{\varphi(\bar{x},\bar{b})\}, q_0)=n^\ast
\]
and by ~\ref{tens:bstar},
\[
\mbox{rk}{}^1_{\varphi(\bar{x},\bar{y})}
(p_0,
q_0\cup\{(\neg\exists\bar{x})(\varphi(\bar{x},\bar{y})\wedge
\varphi(\bar{x},\bar{b})\})<n^\ast.
\]
We have that
$(\neg\exists\bar{x})[\varphi(\bar{x},\bar{y})\wedge
\varphi(\bar{x},\bar{b}]\in q$, hence $q_0\cup
\{(\neg\exists\bar{x})[\varphi(\bar{x},\bar{y})\wedge
\varphi(\bar{x},\bar{b})]\}\subseteq q$, in contradiction with monotonicity
and $\mbox{rk}{}^1_{\varphi(\bar{x},\bar{y})} (p\rest
M_1,q)=n^\ast$.
\end{prf}

\section{More on \sop{2}, \sop{3} and $\initial^*_{\lambda}$-order}
We try to find a connection between the syntactic
properties \sop{2},\sop{3} and the semantic property of
being $\initial^*_{\lambda}$-maximal. Our guess is that
$\initial^*_{\lambda}$-maximality should be equivalent to
one of the above order properties (maybe both), but all we
manage to prove here is \sop{3} \then
$\initial^*_{\lambda}$-maximality \then \sop{2}. We also
give a weak ``local'' result in the other direction.

First we generalize the definitions from \scite{DjSh692},
of $\initial^*_{\lambda}$-maximality, making them local as
well as global.

\begin{dfn}
\begin{enumerate}
\item
    For given (complete first order theories) $T_1,T_2$ and
cardinals $\lambda \ge \mu >\kappa,\mu \ge
\theta > |T_1| + |T_2| + \aleph_0$
\mr
\item [$(a)$]  $T_1 \initial^*_{< \lambda,<\mu,\kappa, < \theta} T_2$
means that there is a (complete first order theory) $T^*$
and interpretations $\bar \varphi_1,\bar \varphi_2$ of
$T_1,T_2$ in $T^*$ respectively, $|T^*| < \theta$ such
that: {\roster \item $\boxtimes^{< \lambda,<
\mu,\kappa}_{T_*,\bar
\varphi_1,\bar \varphi_2}$   if $M$ is a $\kappa$-saturated model of
$T^*$ and $M_\ell = M^{[\bar \varphi^\ell]}$ for $\ell =
1,2$ and $M_2$ is $\lambda$-saturated (model of $T_2$),
\ub{then} $M_1$ is $\mu$-saturated
\eendroster}
\item [$(b)$]  $(T_1,\vartheta_1(\bar x,\bar y)) \initial^*_{< \lambda,<
\mu,< \kappa} (T_1,\vartheta_2(\bar x,\bar y))$ means that
$\vartheta_\ell(\bar x,\bar y) \in L(\tau_{T_\ell})$ and
that there is a $T^*$ and interpretations $\bar
\varphi_1,\bar \varphi_2$ of $T_1,T_2$ in $T^*$
respectively, $|T^*| < \mu$ such that $\boxtimes^{<
\lambda,<
\mu,\kappa}_{T^*,\vartheta_1,\vartheta_2,\bar \varphi^1,\bar \varphi^2}$ if
$M$ is a $\kappa$-saturated model of $T^*$ and $M_\ell =
M^{[\bar
\varphi^\ell]}$ for $\ell =1,2$ and $M_2$ is $(\lambda,\vartheta_1(\bar
x,\bar y))$-saturated (see 3 below), \ub{then} $M_i$ is
$(\mu,\vartheta_2)$-saturated.
\ermn
\item
    Instead ``$< \lambda^+$" we may write ``$\lambda$", and
instead ``$< \mu^+$" we may write $\mu$, instead $<
\theta^+$ we may write $\theta$.  If we omit $\mu$ we mean
$\mu = \lambda$, and if we write $\kappa=0$ then
``$\kappa$-saturated" becomes the empty demand, if we omit
$\theta$ we mean $|T_1| + |T_2| + \aleph_0$ and if we omit
$\kappa$ and $\theta$ then we mean that $\mu =
\lambda,\theta = |T_1| + |T_2| + \aleph_0$.
\item
    We say
$M$ is $(\lambda,\Delta)$-saturated when: if $p \subseteq
\{\vartheta(\bar x;\bar a):\vartheta(\bar x;\bar y) \in \Delta, \bar a \in
{}^{\ell g(\bar y)}M\}$ is finitely satisfiable of
cardinality $<
\lambda$ then $p$ is realized in $M$.
If $\Delta = \{\vartheta(\bar x,\bar y)\}$ we may write
$\vartheta(\bar x,\bar y)$ instead of $\Delta$.
\item
    If $T_1,T_2$ are not necessarily complete, then above $T^*$ is
not necessarily complete and we demand: if $M_1 \models
T_1,M_2 \models T_2$ then there is $M \models T^*$ such
that $M^{[\bar \varphi_\ell]}
\models Th(M_\ell)$ for $\ell=1,2$.

\item
    We say $T$ is $\initial^*_{\lambda,\kappa}$-maximal if
$|T'| <
\lambda \Rightarrow T' \initial^*_{\lambda,\kappa} T$.  We say
$(T,\vartheta(\bar x;\bar y))$ is
$\initial^*_{\lambda,\kappa}$-maximal if $|T'| < \lambda
\and \vartheta'(\bar x';\bar y')) \in L(\tau_{T'})
\Rightarrow (T',\vartheta'(\bar x';\bar y'))
\initial^*_{\lambda,\kappa} (T,\vartheta(\bar x;\bar y))$.

\end{enumerate}
\end{dfn}

\definition
\begin{enumerate}
\item
    $T_{tr}$ is the theory of trees
(i.e. the vocabulary is $\{<\}$ and the axioms state that
$<$ is a partial order and $\{y:y < x\}$ is a linear order
for every $x$), so $T_{tr}$ is not complete, and let
$\vartheta_{tr}(x,y) = (y < x)$.
\item
    $T^*_{tr}$ is the
model completion of $T_{tr}$.
\item
    $T_{\text{ord}}$ is the
theory of linear orders, $T^*_{\text{ord}}$ is its model
completion (i.e. the theory of dense linear order without
endpoints).
\end{enumerate}
\eenddefinition

We note connection to previous works and obvious properties
\proclaim
\begin{enumerate}
\item
    $T_1 \initial^*_{\lambda,\mu,0} T_2$ is $T_1
\initial^*_{\lambda,\mu} T_2$ of
\cite{DjSh692}.
\item
    $T_1 \initial^*_{\lambda,\lambda;< \kappa} T_2$ implies
$T_1
\initial^*_{\lambda,\kappa} T_2$ of \cite[2.x,p.xxx]{Sh500}.
\item
    $\initial^*_{\lambda,\mu;\kappa,\theta}$ has the obvious
monotonicity properties: if $T_1 \initial^*_{<
\lambda_1,<\mu'_1;<
\kappa_1,< \theta_1} T_2$ and $\lambda_2 \ge \lambda_1,\mu_2 \le
\mu_1,\kappa_2 \ge \kappa_1,\theta_2 \ge \theta_1$ then $T_1
\initial^*_{< \lambda_2,< \mu_2;< \kappa_2,< \theta_2} T_2$.
\item
    $T \initial^*_{\lambda,\mu;\kappa,\theta} T$ if $|T| <
\theta,\lambda \ge \mu > \kappa,\mu \ge \theta$.
\item
    If $\mu$ is a limit cardinal, then $T_1
\initial^*_{< \lambda,< \mu;< \kappa,< \theta} T_2$ iff for every
$\mu_1 < \mu,\mu_1 \ge \kappa$ we have

$$
T_1 \initial^*_{< \lambda,< \mu_1;< \kappa,< \theta} T_2.
$$
\item
    Similar results hold for $(T_\ell,\vartheta_\ell(\bar
x;\bar y))$.
\end{enumerate}
\eendproclaim

\begin{prf}  Easy.
\end{prf}

\proclaim
\begin{enumerate}
\item
    Assume $T_1
\initial^*_{< \lambda,< \mu;< \kappa,< \theta} T_2$.  \ub{Then}
for any theory $T^*$, we can find $T^{**} \supseteq T^*$
complete $|T^{**}| < (|T^*|^{|\tau(T_1)|+|\tau(T_2)|})^+ +
\theta$ such that: for any interpretations $\bar
\varphi_1,\bar \varphi_2$ of $T_1,T_2$ in $T^{**}$
respectively the definition  of $T_1
\initial^*_{< \lambda,< \mu;< \kappa,< \theta} T_2$ holds.
\item Assume $\tau(T_1),\tau(T_2)$ are disjoint.  Then $T_1
\initial^*_{< \lambda,< \mu;< \kappa,< \theta} T_2$ if for any $T
\supseteq T_1 \cup T_2$ there is $T^* \supseteq T$ as demanded in
Definition for the trivial interpretations
$M^{[\bar
\varphi^\ell]}$ is the $\tau(T_\ell)$-reduct.
\end{enumerate}
\eendproclaim

\begin{prf} Easy.
\end{prf}

Now we will show that $T^*_{tr}$ is
$\initial^*_\lam$-maximal for every \lam \; big enough, and
conclude that \sop{3} \then $\initial^*_\lam$-maximality.
The last result appears already in \cite{Sh500}, theorem
(2.9), but the proof is not full - in fact, the proof shows
the following theorem:

\begin{thm}
\label{thm:sopthree}
    Any theory $T$, $\card{T} < \lam$,
    with \sop{3} is $\initial^*_\lam$-above
    $T^*_{ord}$.
\end{thm}
\begin{prf}
    See \cite{Sh500}, (2.12).
\end{prf}

Here we fill the missing part, proving explicitly that
$T^*_{tr}$, and therefore $T^*_{ord}$ are maximal.

\begin{thm}
$T^*_{tr}$ is $\initial^*_\lambda$-maximal for any $\lambda
> \aleph_0$; the witness $T^*$ does not depend on
$\lambda$.
\end{thm}

\begin{rmk}  This continues \cite[Ch.VI,3.x]{Shc}.
\end{rmk}

\begin{prf}  Let $T$ be any complete theory, $|T| < \lambda$ and
$M_1$ a model of $T$.

Let $\Phi = \{\varphi(x,\bar a):\varphi(x,\bar y) \in
L(\tau_T),\bar a \in {}^{\ell g(\bar y)}(M_1)\}$, so
$|\Phi| = \|M_1\|$.  So $M = ({}^{\omega >} \Phi,\initial)$
is a model of $T_{tr}$ and there is a model $M_2$ of
$T^*_{tr}$ of cardinality $\|M_1\|$ extending $M$ such that
every member of $M_2$ is below some member of $M$.

Let $\chi$ be large enough such that $M_1,M_2 \in {\Cal
H}(\chi)$ and we define ${\Cal B}^*$ expanding $({\Cal
H}(\chi),\in)$ by $P_1 = |M_1|,P_2 = |M_2|,P = |M|,Q_0 =
\Phi,<_1 = <^{M_2},< = <_1 \restriction P, m = $ a constant
symbol for a set $M_1, R^{{\Cal B}^*} = R^{M_1}$ for $R
\in \tau_T$ (wlog $\tau(T)$ does not contain any other
predicate mentioned here)

$$
Q = \{(\langle \varphi_\ell(x,\bar a_\ell):\ell < n
\rangle:M_1
\models \exists x [\wedge \varphi_\ell(x,\bar a_\ell)]\}.
$$

$H$ is a partial unary function with domain $Q$ and range
$P_1$, $H(\langle \varphi_\ell(x,\bar a_\ell):\ell < n
\rangle)$ satisfies $\{\varphi_\ell(x,\bar a_\ell):\ell <
n\}$, i.e. ${\Cal B}^*$ satisfies the formula $``m \models
(\exists x)
\bigwedge_{\ell < n}\varphi_{\ell}(x, \bar a_{\ell})''$.

Let $T^* = Th({\Cal B}^*)$, let $\bar \varphi_1$ be the
trivial interpretation of $T$ in $T^*$ (the restriction +
reduct) and $\bar
\varphi_2 = \langle P_2(x),x_0 <_1 x_1 \rangle$
is an interpretation of $T^*_{tr}$. So $T^*,\bar
\varphi_1,\bar \varphi_2$ does not depend on $\lambda$.

Now we assume ${\Cal B}$ is a model of $T^*,N_1 = {\Cal
B}^{[\bar \varphi_1]},N_2 = {\Cal B}^{[\bar \varphi_2]},N_3
= (P^{\Cal B},<^{\Cal B})$ and we aim to show that $(i)$
below implies $(iii)$. We will first show that $(i)
\Rightarrow (ii)$ and use this fact in the proof.
\roster
\itemitem[ $(i)$ ]  $N_2$ is $\lambda$-saturated
\itemitem[ $(ii)$ ]  in $N_3$ every branch has cofinality
$\ge \lambda$,
equivalently: every increasing sequence of length $<
\lambda$ has an upper bound
\itemitem[ $(iii)$ ]  $N_1$ is $\lambda$-saturated.
\eendroster

\ub{Why $(i) \Rightarrow (ii)$?}  If $\langle a_i:i < \delta \rangle$
is $<^{N_3}$-increasing, $\delta < \lambda$ then it is
$<^{N_2}$-increasing hence has a $<^{N_2}$-upper bound $a$
but $(\forall x \in P_2)(\exists y)(x <_1 y \and P(y))$
belongs to $T^*$ so there is $b,a <^{N_2} b \in P^N = N_3$
so $b$ is as required.

So we can assume clause (i) and we shall prove (iii).

\bigskip
Before we proceed, let us note several trivial but
important properties of ${\Cal B}$.
\mr
\item [$(a)$] We can talk inside $\CB$ about a set being a model,
(standard coding of) a formula, a proof, etc. In
particular, we can speak about $m$ (as a model) satisfying
or not satisfying certain sentences. Also, given a formula
with free variables we can speak about substitution of
other variables or parameters into the formula. Given $s
\in \CB$ which is a formula with free variables $\bar x$,
we will allow ourselves to write $s = s(\bar x)$, and if
$\CB$ thinks that substitution of $\bar a \in P_1$ into $s$
will turn it into a true sentence in $m$ as a model, we
will write $m \models s(\bar a)$ or just $s(\bar a)$.

\item [$(b)$] ${\Cal B} \models
\forall z Q_0(z) \ifff ``z$ is a formula with one free variable with
parameters from $P_1''$.
Moreover, suppose $\varphi(x, \bar a)$ is a formula in
$L(\tau_{T})$ s.t. $\bar a \in P_1^\CB$. ${\Cal B^*}$ and
therefore $\CB$ satisfy $(\forall {\bar y} \in P_1)(
\exists ! s \in Q_0)$ such that $(\forall x \in P_1)
\varphi(x, \bar y) \ifff ``m \models s(x, \bar y)''$.
Let us denote by $\ulcorner \varphi(x, \bar a)
\urcorner$ this ``canonical encoding'' of $\varphi(x, \bar a)$ in $Q_0^\CB$.

\item [$(c)$] $\CB \models \forall s P(s) \ifff ``s$ is a finite sequence
of members of $Q_0$, i.e. $(\exists n \in \omega) (s:n \to
Q_0)''$.

\item [$(d)$] For simplicity of notation, given $s \in P^\CB$, we will
write $``z \in s''$ instead of $``z \in Im(s^\CB)''$.

\item [$(e)$] For $z \in P^\CB$, $c \in P_1^\CB$, we write
$z(c)$ meaning $(\forall s \in z) s(c)$.

\item [$(f)$] For every $\varphi(x, \bar a) \in L(\tau_T)$ for
$a \in P_1^\CB$, there exists an element of $P^\CB$
corresponding to the finite sequence $\langle \varphi(x,
\bar a) \rangle$. We denote this element by $\langle
\ulcorner \varphi(x, \bar a) \urcorner \rangle$. Moreover,
$\CB \models \exists x (P_1(x) \land \varphi(x, \bar a))
\to Q(\langle \ulcorner \varphi(x, \bar a) \urcorner
\rangle)$.

\ermn

\begin{subclm}
\begin{enumerate}
\item Suppose $\CB \models Q(z)$. Then $\CB \models \forall w
(Q(w) \land z<w) \to z(H(w))$.

\item Let $\varphi(x, \bar a)
\in L(\tau_T)$ and suppose $\CB
\models \exists x P_1(x) \land \varphi(x, \bar a)$. Then
$\CB \models \forall z (Q(z) \land \langle \ulcorner
\varphi(x, \bar a) \urcorner \rangle < z)
\to \varphi(H(z), \bar a)$.

\end{enumerate}
\end{subclm}

\begin{prf}
\begin{enumerate}
\item Trivial as ${\Cal B^*}$ satisfies it.

\item Let $z^* =
\langle \ulcorner \varphi(x, \bar a)
\urcorner \rangle$. First, $Q(z^*)$ holds by $f$ above. By
(1), $z^*(H(z))$ holds for each $z \in Q^\CB, z^* < z$. Now
by $b$ and $f$ above, $\CB \models \forall x P_1(x) \to
(z^*(x) \ifff
\varphi(x, \bar a))$. As $\CB \models Range(H) \subseteq P_1$, we are done.
\end{enumerate}
\end{prf}

We now proceed with the proof (i) $\Longrightarrow$ (iii).
So let $p$ be a 1-type in $N_1$ of cardinality $< \lambda$,
so let $p = \{\varphi_\beta(x,\bar a_\beta):\beta <
\alpha\}$ with $\alpha <
\lambda$, $\bar a_\beta \in N_1 \forall \beta$.
Without loss of generality $p$ is closed under conjunction,
i.e. for every $\varepsilon,\zeta < \alpha$ for some $\xi <
\alpha$ we have $\varphi_\xi(x,\bar a_\xi) =
\varphi_\varepsilon(x,\bar a_\varepsilon) \wedge
\varphi_\zeta(x,\bar a_\zeta)$.  We shall now choose by
induction on $\beta \le \alpha$ an element $b_\beta$ of $N$
such that
\mr
\item [$(A)$]  $b_\beta \in P^{\Cal B} = N_3$ moreover $b_\beta
\in Q^{\Cal B}$ and $\gamma < \beta \Rightarrow b_\gamma <^{N_3}
b_\beta$

\item [$(B)$]  if $\gamma < \beta$ then
${\Cal B} \models (\forall z)(Q(z) \wedge (b_\beta \le z)
\rightarrow \varphi_\gamma(H(z),\bar a_\gamma))$

\item [$(C)$]  if $\gamma < \alpha$ (but not necessarily $\gamma <
\beta$) then ${\Cal B} \models (\exists z)[Q(z) \wedge (b_\beta
\le z) \wedge (\forall y)(Q(y) \wedge z \le y \rightarrow
\varphi_{\gamma}(H(y),\bar a_\gamma))]$.
\ermn

If we succeed then $H^{\Cal B}(b_\alpha)$ is as required.

\ub{Case 1}:  $\beta = 0$.

Define $b_0 = \langle \rangle$ (the element of $P^\CB$
corresponding to the empty sequence). Clearly $\CB \models
Q(b_0)$, i.e.  the demand $(A)$ holds. $(B)$ holds
trivially. Why does $(C)$ hold? Let $\gamma < \alpha$. $\CB
\models \exists x \varphi_\gamma(x, \bar a_\gamma)$
therefore denoting $z^*_\gamma = \langle \ulcorner
\varphi_\gamma(x, \bar a_\gamma)
\urcorner \rangle$, we have $\CB \models Q(z^*_\gamma) \land b_0 < z^*_\gamma$.
Now we finish by part (2) of the subclaim.

\ub{Case 2}:  $\beta = \upsilon + 1$.

$\CB$ satisfies the sentence saying that for every $\eta
\in Q$ and $\bar y \in P_1$ there exists an element of $P$
that we denote by $Conc_\upsilon(\eta, \bar y)$
corresponding to $\eta \char 94 \langle \ulcorner
\varphi_\upsilon(x,\bar y) \urcorner \rangle$. We define
$b_\beta = Conc_\upsilon(b_\upsilon, \bar a_\upsilon)$. Now
we have to check $(A)$ - $(C)$.

\mr
\item [$(A)$] By the induction hypothesis, clause $(C)$ holds for $b_\upsilon$
and $\upsilon$ (standing for $b_\beta$ and $\gamma$ there).
Therefore $\CB \models \exists z \in Q (b_\upsilon \le z)
\land
\varphi_\upsilon(H(z),\bar a_\upsilon))$. But ${\Cal B^*}$ (and so $\CB$) satisfies
that $\forall \bar y \in P_1$ if there exists $z \in Q$
s.t. $\varphi_\upsilon(H(z), \bar y)$ holds, then
$Conc_\upsilon(z, \bar y)$ is an element of $Q$ (as in
${\Cal B^*}$ the assumption means that there exists an
element of $m$ satisfying all the formulae in $z$ plus
$\varphi_\upsilon(x, \bar y)$). So we get the required.

\item [$(B)$] is clear as by the induction hypothesis,
$\varphi_\zeta(H(z),\bar a_\zeta)$ holds for every $\zeta <
\upsilon$, $b_\beta \le z$ (recall that $b_\upsilon \le
b_\beta$). As for $\varphi_\upsilon(x, \bar a_\upsilon)$,
${\Cal B^*}$ clearly satisfies that for every $z \in Q,
\bar y \in P_1$, if $b = Conc_\upsilon(z, \bar y)$ is in
$Q$ then $\varphi_\upsilon(H(z), \bar y)$ holds $\forall z
\in Q, b \le z$.

\item [$(C)$] Let $\zeta < \alpha$. As $p$ is closed under conjunctions,
for some $\xi$, $\varphi_\gamma(x,\bar a_\gamma) \land
\varphi_\zeta(x,\bar a_\zeta) = \varphi_\xi(x,\bar a_\xi)$.
Now we apply clause $(C)$ holding for $b_\upsilon$ to
$\gamma = \xi$ and get $z \in Q, b_\upsilon \le z$ with
$H(z)$ satisfying both $\varphi_\upsilon(x,\bar
a_\upsilon)$ and $\varphi_\zeta(x,\bar a_\zeta)$. Once
again using the satisfaction by ${\Cal B}$ of natural
sentences, we show that $b = Conc_\zeta(b_\beta, \bar
a_\zeta)$ is in $Q$, $b_\beta \le b$ and $\forall z \in Q$
which is above $b$, $\varphi_\zeta(x, \bar a_\zeta)$ holds,
i.e. b is as required.

\ermn

\ub{Case 3}:  $\beta = \delta$ limit.

By our present assumption, clause (i), and therefore clause
(ii), hold. Hence there is $b \in P^\CB$ which is an upper
bound to $\{b_\gamma:\gamma < \beta\}$.  Now $\CB$
satisfies ``for every element $z$ of $P$ there is a $y \le
z$ which is in $Q$ and $x \le z \and Q(x) \to x \le y$".
Apply this to $b$ for $z$ and get $b'_\delta$ for $y$.  So
$b'_\delta \in Q$ and $\gamma < \delta \Rightarrow b_\gamma
\le b'_\delta$, as required in clauses $(A)$ +$(B)$ but not
necessarily $(C)$.

Define for each $\zeta < \alpha$ a formula $\psi_\zeta(w,
\bar a_\zeta) = (\exists z)(w \le z \land Q(z) \wedge
(\forall y) (z \le y \land Q(y) \rightarrow
\varphi_\zeta(H(y),\bar a_\zeta))$ Now we find $c_\zeta$
(for $\zeta < \alpha$) such that:

\mr
\item [$(a)$]   $c_\zeta \in Q^{\Cal B},c_\zeta \le b$

\item [$(b)$]  $\psi_\zeta(c_\zeta, a_\zeta)$ holds.

\item [$(c)$]  under $(a)$ + $(b)$, the element $c_\zeta$ is maximal.
\ermn

Why do $c_\zeta$ exist? $\CB$ satisfies ``for every element
$s$ of $P$ there is a $w \le s$ which satisfies
$\psi_\zeta(w, \bar a_\zeta)$, is in $Q$ and $(x \le s
\land \psi_\zeta(x, \bar a_\zeta) \land Q(x))
\to (x \le w)$".

By the induction hypothesis we have:

$$
\gamma < \delta, \zeta < \alpha \Rightarrow b_\gamma <^{N_3} c_\zeta.
$$

Clearly it suffices to find $b_\delta$ satisfying
$Q(b_\delta)$ and $b_\gamma <^{N_3} b_\delta <^{N_3}
c_\zeta$ for $\gamma < \delta,\zeta < \alpha$. As $N_3
\restriction \{c:c \le b\}$ is linearly ordered, this
follows from $N_2$ being $\lambda$-saturated.
\end{prf}

\proclaim
\label{prp:ord}
\begin{enumerate}

\item For every $T^*$, there is $T^{**} \supseteq T^*,|T^{**}|
= |T^*| +
\aleph_0$ such that for every model ${\Cal B}$ of $T^{**}$ we have
\mr
\item [$(a)$]  for any $\lambda$, the following are equivalent
{\roster
\itemitem [$(\alpha)$ ]  if $\bar \varphi_1$ is an interpretation of
$T^*_{tr}$ in ${\Cal B}$ (possibly with parameters) then
${\Cal B}^{[\bar \varphi_1]}$ is $\lambda_{tr}$-saturated

\itemitem[ $(\beta)$ ]  if $\bar \varphi_2$ is an interpretation of
$T_{\text{ord}}$ in ${\Cal B}$ (possibly with parameters)
then ${\Cal B}^{[\bar \varphi_2]}$ is $\lambda$-saturated
\eendroster}

\item [$(b)$]  for any $\lambda$, the following are equivalent
{\roster
\itemitem[ $(\alpha)$ ]  if $\bar \varphi_1$ is an interpretation of
$T_{tr}$ in ${\Cal B}$ (possibly with parameters) then in
${\Cal B}^{[\bar \varphi_1]}$, every branch with no last
element has cofinality $\ge \lambda$

\itemitem[ $(\beta)$ ]  if $\bar \varphi^*_2$ is an interpretation of
$T_{\text{ord}}$ in ${\Cal B}$ (possibly with parameters)
then in ${\Cal B}^{[\bar \varphi_2]}$ there is no Dedekind
cut $(I_1,I_2)$ with both cofinalities $< \lambda$ and at
least one $\ge \aleph_0$.
\eendroster}
\eendroster
\end{enumerate}
\eendproclaim

\begin{prf}  Easy.
\end{prf}

\begin{cor}\label{cor:max}
\begin{enumerate}
\item $T^*_{\text{ord}}$ is
$\initial^*_{\lambda}$-maximal.

\item If $|T| <
\lambda$ and $T$ has SOP$_3$ then $T$ is
$\initial^*_{\lambda}$-maximal.
\end{enumerate}
\end{cor}

\begin{prf}
\begin{enumerate}
\item Follows from ~\ref{prp:ord}
\item By (1) and ~\ref{thm:sopthree}.
\end{enumerate}
\end{prf}

\begin{qst}
    Is the other direction of ~\ref{cor:max} (2) true?
\end{qst}

\begin{rmk}
    We present later a proof of a weaker version of the
    other direction: we get \sop{2} instead of \sop{3}.
\end{rmk}

We would like to prove a result similar to
~\ref{thm:sopthree} for \sop{2} (or to show maximality in
some other way), but unfortunately right now we only can
present the following local theorem:

\begin{thm}  If $T$ has SOP$_2$ as exemplified by
$\vartheta(\bar x;\bar y)$, \ub{then}
$(T^*_{tr},\vartheta_{tr}(x;y))
\initial^*_\lambda (T,\vartheta(\bar x;\bar y))$ for any $\lambda \ge
|T|+\aleph_0$ regular.
\end{thm}

\begin{prf}  We can find a model $M_1$ of $T^*_{tr}$ and model $M_2$ of
$T$ and $\bar a_b \in {}^{\ell g(\bar y)}M_2$ for $b \in
M_1$ such that:
\mr
\item [$(\alpha)$]  if $M_1 \models b_0 < \ldots < b_{n-1}$ then
$\{\vartheta(\bar x,\bar a_{b_\ell}):\ell < n\}$ is
satisfiable in $M_2$
\item [$(\beta)$]  if $b_1,b_2$ are incomparable in $M_1$ then
$$
M_2 \models \neg(\exists \bar x)(\vartheta(\bar x,\bar
a_{b_1}) \and
\vartheta(\bar x,\bar a_{b_2}))
$$
\item [$(\gamma)$] for no $\bar d \in {}^{\ell g(\bar x)}(M_2)$ is
$\{b \in M_1:M_2 \models \vartheta(\bar d,\bar a_b)\}$
unbounded in $M_1$ (note that by $(\beta)$ it is always
linearly ordered in $M_1$, therefore $(\gamma)$ means that
for each $\bar d \in {}^{\ell g(\bar x)}(M_2)$, there
exists an element of $M_1$ which is above every $b$
satisfying $\vartheta(\bar d,\bar a_b)$).

[How? Choose by induction on $n,(M_{1,n},M_{2,n},\langle
\bar a_b:b \in M_{1,n} \rangle : n < \omega)$ such that:
\roster
\itemitem[ $(a)$ ]   $M_{1,n}$ is a model of $T^*_{tr}$
\itemitem[ $(b)$ ]   $M_{2,n}$ is a model of $T$
\itemitem[ $(c)$ ]   $M_{1,n} \prec M_{1,n+1}$ moreover, every branch of
$M_{1,n}$ has an upper bound in $M_{1,n+1}$
\itemitem[ $(d)$ ]   $M_{2,n} \prec M_{2,n+1}$
\itemitem[ $(e)$ ]   $\bar a_b \in {}^{\ell g(\bar y)}(M_{2,n})$ for
$b \in M_{1,n}$
\itemitem[ $(f)$ ]  clauses $(\alpha),(\beta)$ hold
\itemitem[ $(g)$ ]   if $b \in M_{1,n+1}$ and $[b' \in M_{1,n}
\Rightarrow M_{1,n+1} \models \neg(b < b')]$ then $\vartheta(\bar x,\bar
a_b)$ is not satisfied by any sequence from $M_{1,n}$.
\eendroster
There is no problem to carry the definition.

Now $M_1 =
\dbcu_n M_{1,n},M_2 = \dbcu_n M_{2,n}$ and $\langle \bar
a_b:b \in M_1 \rangle$ are as required above.]
\ermn
Now let $\chi$ be such that $M_1,M_2 \in {\Cal H}(\chi)$,
\wilog \,\, $\tau_T = \tau(M_2)$, $\{<\} = \tau(T_{tr})
= \tau(M_1)$ and $\{\in\}$ are pairwise disjoint.  Now we
define a model ${\Cal B}_0$.

Its universe is ${\Cal H}(\chi)$
\ub{relation} $\in$ (membership)

$P_1 = |M_1|$,

$P_2 = |M_2|$

$R = R^{M_\ell} \text{ if } R \in \tau(M_\ell),\ell \in
\{1,2\}$
$F_\ell$ (for $\ell < \ell g(\bar y)$) a partial unary
function such that: $b \in M_1 \Rightarrow \langle
F_\ell(b):\ell < \ell g(\bar y)
\rangle = \bar a_b$.

Let $T^* = Th({\Cal B}_0)$.  For the obvious $\bar
\varphi$ and $\bar \psi,T^*$ is $(T,T_{tr})$-superior and $|T^*| = |T| +
\aleph_0$.  Assume $\lambda = \text{ cf}(\lambda) > |T^*|$.

So let ${\Cal B}$ be a model of $T^*$ such that $M'_2 =
{\Cal B}^{[\bar \varphi]}$, the model of $T$ interpreted in
it, is $\lambda^+$-saturated. It will be enough to prove
that $M'_1 = {\Cal B}^{[\bar
\psi]}$ satisfies: for every branch of cofinality $\theta \le \lambda$
there exists an upper bound. So let $\{b_i:i < \theta\}$ be
$<^{M_1}$-increasing let $\bar c_i =
\langle F^{\Cal B}_\ell(b_\ell):\ell < \ell g(\bar y) \rangle$.  Hence
for any $n < \omega,i_0 < \ldots < i_{n-1} < \theta$ we
have $M'_2
\models (\exists \bar x)[\dsize \bigwedge_{m<n} \vartheta(\bar x, \bar
c_i)]$ because ${\Cal B}_0 \models (\forall
z_0,\dotsc,z_{n-1})[\dsize
\bigwedge_{k < n} P_1(z_k) \Rightarrow z_0 < z_1 < \ldots <
z_{n-1} \rightarrow (\exists \bar x) \dsize \bigwedge_{m<n}
\vartheta(\bar x,\langle F_\ell(z_m):\ell < \ell g(\bar y) \rangle)]$.

So $\{\vartheta(\bar x,\bar c_i):i < \theta\}$ is finitely
satisfiable in $M'_2$ hence some $\bar d \in {}^{\ell
g(\bar x)}(M'_2)$ realizes it.  Now we claim that $\{b \in
M'_1:{\Cal B} \models \vartheta(\bar d,\bar a_b)\}$ is
bounded in $M'_1$. Why? Recall that by clause ($\gamma$)
${\Cal B_0}$ satisfies: for every $\bar x \in {}^{\ell
g(y)} P_2$ there exists $z \in P_1$ such that $z$ is
${<}^{\Cal B}$ - above all the elements $w \in P_1$
satisfying $\vartheta(\bar x, \bar a_w)$. Therefore ${\Cal
B}$ satisfies this sentence, and applying it to $\bar d \in
{}^{\ell g(\bar x)}(M'_2)$, we get $b^* \in M'_1$
 - the required bound. As for each $i < \theta$,
$\vartheta(\bar d, \bar{a_{b_i}})$ holds, clearly ${\Cal B}
\models b_i < b^*$ for all $i$, and we are done.
\end{prf}

The next goal is to complete the proof started in
\cite{DjSh692} of the fact that $\initial^*$-maximality
implies \sop{2}. In \cite{DjSh692} a property was defined -
$\initial^{**}_\lam$ - maximality, which is closely related
to $\initial^{*}_\lam$ - maximality and it was shown in
theorem (3.4) that every $T$ which is $\initial^{**}_\lam$
- maximal for some (every) big enough regular \lam, has an
order property similar to \sop{2}, that we call $SOP''_2$
(see ~\ref{dfn:nsop21}). We answer the question (3.8)(3)
from [DjSh692] showing that $SOP''_2$ is equivalent to
\sop{2} (for a theory).

So assuming that $T$ is $\initial^{*}_{\lam^+}$ - maximal
for some regular
\lam \; satisfying $2^\lam =
\lam^+$, we get by \scite{DjSh692}, claim (3.2), $T$ is
$\initial^{**}_\lam$ - maximal, so it has $SOP''_2$, and
therefore \sop{2}.

\begin{thm}
Let $T$ be a theory.
\begin{enumerate}
\item  Suppose $\vartheta(\bar x,\bar y)$ exemplifies SOP$_2$ in
$T$. Then $\vartheta(\bar x;\bar y)$ exemplifies SOP$''_2$
in $T$ as well.

\item  Suppose $\vartheta(\bar x,\bar y)$ exemplifies SOP$''_2$ in
$T$. Then for some $k$, $\vartheta^{<k>}(\bar x;\bar y)$
exemplifies SOP$_2$ in $T$ (where $\vartheta^{<k>}(\bar
x;\bar y^{<k>}) = \dsize \bigwedge_{\ell < k}
\vartheta(\bar x;\bar y_\ell)$).
\end{enumerate}
\end{thm}

\begin{prf}
\begin{enumerate}
\item is easy.

\item  Denote ${\Cal
I}^n_\lambda = \{\bar \eta:\bar \eta = \langle
\eta_\ell:\ell \le n \rangle,\eta_\ell \initial \eta_{\ell +1}$;
and $\eta_\ell \in {}^{\lambda >}2\}$. So assume
$\vartheta(\bar x;\bar y)$ has SOP$''_2$ as exemplified by
$n,\bar{\mathbf a} = \langle a_{\bar
\eta}:\bar \eta \in {\Cal I}^n_\omega \rangle$.  Without loss of
generality $\langle \bar a_{\bar \eta}:\bar \eta \in {\Cal
I}^n_\omega
\rangle$ is tree indiscernible in the relevant sense:
\dcontz{\et}, \dconto{\et} look the same over \et \,
($2-fbti$ from ~\ref{dfn:indisc}). We can assume this by
~\ref{fct:thinning} (for more details, see \cite{DjSh692},
claim (2.14)).


For $\nu \in {}^{\omega \ge}2$ let $p_\nu =
\{\vartheta(\bar x,\bar a_{\bar \eta}):\bar \eta = \langle
\eta_\ell:\ell < n
\rangle,\eta_\ell < \eta_{\ell +1} \initialeq \nu\}$ so
\mr
\item [$\circledast_1$]  $p_\eta$ for $\eta \in {}^\omega 2$ is
consistent (in ${\mathfrak C}_T$).
\ermn
Let

\[
\begin{array}{ll}
&\Xi = \bigl\{(h,\Upsilon):h \text{ is a one-to-one mapping
from}
    {}^{n\ge}m \text{ to } {}^{\omega >}2 \\
&\text{preserving } \initial,\perp \text{ and } \Upsilon
\subseteq
    {}^n m \text{ and there is } \\
&\langle \nu^*_\eta:\eta \in \Upsilon \rangle,h(\eta)
    \initial
    \nu^*_\eta \in {}^\omega 2 \text{ for } \eta \in{}^n m \\
&\text{ such that } \cup\{p_{\nu^*_\eta}:\eta \in
\Upsilon\}
    \text{ is inconsistent} \bigr\}.
\end{array}
\]

Now
\mr
\item [$\circledast_2$]  $\Xi$ is nonempty

[Why?  By the definition of $SOP''_2$, clause $(b)$, choose
$\Upsilon = ^n m$ ]
\ermn
Choose $(h^*,\Upsilon^*) \in \Xi$ with $|\Upsilon^*|$ of
minimal cardinality and $\langle \nu^{*}_\eta:\eta \in
\Upsilon^* \rangle$ as there.
By $\circledast_1$ clearly $|\Upsilon^*| \ge
2$. So choose $\eta_0 \ne \eta_1$ from $\Upsilon^*$ with
\mbox{$\nu^* = \nu^*_{\eta_0} \cap \nu^*_{\eta_1}$}
\mbox{$( = h(\et_0) \cap h(\et_1))$}
being of maximal \fbox{????} length and let $k^* = \ell
g(\nu^*)$. We can find $\ell^* < \omega$ sufficiently large
such that $\cup
\{p_{\nu^*_\eta
\restriction
\ell^*}:\eta \in \Upsilon^*\}$ is inconsistent.  We choose by induction on $i <
\omega$ for every $\rho \in {}^\ell 2$, a sequence $\nu_\rho \in
{}^{\omega >} 2$ by $\nu_{<>} = \nu^*,\nu_{\rho \char 94
<j>} =
\nu_\rho \char 94 h(\eta_j)$.

Lastly for $\rho \in {}^{\omega >}2 \in \{<>\}$ let
$\vartheta^*(\bar x,\bar b^*_\rho)$ be the conjunction of
\[
\begin{array}{ll}
&\bigcup \bigl\{ p_{\nu^*_\eta \restriction
\ell^*}:\eta
    \in \Upsilon^* \backslash \{\eta_0,\eta_1\} \}\cup
    \{\vartheta(\bar x,\bar a_{\bar
    \eta}):\bar \eta = \langle \eta_\ell:\ell \le n \rangle, \\
&\eta_\ell \initial \eta_{\ell +1} \initialeq \nu_\rho
    \text{ and } (\forall \ell \le n)[\ell g(\eta_\ell) \notin [k,\ell
    g(\nu_\rho) - \ell^*)] \\
&\text{(this condition is empty if } \ell g(\rho) = 1)
    \bigr\}.
\end{array}
\]

Now if $\rho^* \in {}^\omega 2$ then $\{\vartheta^*(x,\bar
b_\rho):\rho \initial \rho^*\}$ is consistent as all its
members are conjunctions of formulas from
$$
\cup \{p_{\nu^*_\eta}:\eta \in \Upsilon^* \backslash \{\eta^*_0,\eta^*_1\}
\cup p_{\rho^*}
$$
and this is consistent as otherwise $(h^* \restriction
(\Upsilon^*
\backslash \{\eta^*_0,\eta^*_1\}) \cup \{\langle
\eta^*_0,\rho^* \rest \l^{**}
\rangle\},\Upsilon^* \backslash \{\eta^*_1\})$ belongs to \;
$\Xi$ \; for some $\l^{**}$, thus contradicting the choice
of $(h^*,\Upsilon^*)$, i.e. with minimal $|\Upsilon^*|$.

Lastly if $\rho_0,\rho_1 \in {}^{\omega >} 2$ are
$\initial$-incomparable then $\{\vartheta^*(\bar x;\bar
b_{\rho_0}),\vartheta^*(\bar x;\bar b_{\rho_1})\}$ is
inconsistent: we know that

\begin{tens_equ}
\label{tens:incon}
\[
\begin{array}{ll}
&\bigcup \{ p_{\nu^*_\eta \restriction
\ell^*}:\eta
    \in \Upsilon^* \backslash \{\eta_0,\eta_1\} \}\cup
    \{\vartheta(\bar x,\bar a_{\bar
    \eta}):\bar \eta = \langle \eta_\ell:\ell \le n \rangle, \\
&\eta_\ell \initial \eta_{\ell +1} \initialeq
    \nu^*_{\et_0} \rest \l^* \} \cup \{\vartheta(\bar x,\bar a_{\bar \eta}):
    \bar \eta = \langle \eta_\ell:\ell \le n \rangle,
    \eta_\ell \initial \eta_{\ell +1} \initialeq
    \nu^*_{\et_1} \rest \l^*
    \}.
\end{array}
\]
\end{tens_equ}

is inconsistent (by the choice of $(h^*,\Upsilon^*) \in
\Xi$\; and the choice of $\l^*$). Now, by the fact that
\mbox{$\nu^* = \nu^*_{\eta_0} \cap \nu^*_{\eta_1}$} was
chosen to be maximal among other pairs in $\Upsilon^*$, we
see that if
\[
    \bar \eta_0 = \langle \eta^0_{\ell} : \ell
    \le n \rangle, \text{where for each \;} \l, \;
    \eta^0_{\ell} \initial \eta^0_{\ell +1} \initialeq
    \nu^*_{\et^*_0} \rest \l^*
\]
and
\[
    \bar \eta_1 = \langle \eta^1_{\ell} : \ell
    \le n \rangle, \text{where for each \;} \l, \;
    \eta^1_{\ell} \initial \eta^1_{\ell +1} \initialeq
    \nu^*_{\et^*_1} \rest \l^*
\]
while
\[
    \bar \eta_3 = \langle \eta^3_{\ell} : \ell
    \le n \rangle, \text{where for each \;} \l, \;
    \eta^2_{\ell} \initial \nu^*_{\et^*} \text{for some} \;
    \eta^* \in \Upsilon^* \setminus \set{\eta^*_0,\eta^*_1}
\]
then
\begin{tens_equ}
\label{tens:sim}
\[
    \bar \eta_1 \cont \bar \eta_2 \cont \bar \eta_3 \equiv
    \bar \varsigma_1 \cont \bar \varsigma_2 \cont \bar
    \eta_3
\]
\end{tens_equ}
where $\bar \varsigma_j = \langle \varsigma^j_{\ell} : \ell
\le n \rangle$ and
\[
\begin{array}{ll}
    & \varsigma^j_{\ell} = \eta^j_{\l},\; \tif \lg(\eta^j_{\l})
        \le k^* \\
    & \varsigma^j_{\ell} = \nu_{\rho_j} \rest
        [\lg(\nu_{\rho_j}) -
        (\l^*-\lg(\eta^j_{\l}))],\;
        \text{otherwise}
\end{array}
\]
In simpler words: we replace every $\eta^j_\l$ (an initial
segment of $\nu_{\eta_j} \rest \l^*$) whose length is
bigger than $k^*$ (in particular, it is not below any
element in the image of $\Upsilon^*$ other than
$\nu_{\eta_j}$ itself ) by an appropriate initial segment
of $\nu_{\rho_j}$, and get a similar sequence over the
image of $\Upsilon^* \setminus \set{\eta^*_0, \eta^*_1}$.

Now, by indiscernibility of \seq{\a{\bar eta}}, the
definition of $\thet^*(\x,\bar b^*_\rho)$,
~\ref{tens:incon} and ~\ref{tens:sim}, we conclude
$\{\vartheta^*(\bar x;\bar b_{\rho_0}),\vartheta^*(\bar
x;\bar b_{\rho_1})\}$ is also inconsistent.

\end{enumerate}
\end{prf}

\end{document}